# Normalized solutions to fractional mass supercritical NLS systems with Sobolev critical nonlinearities

Jiabin Zuo[1,2] · Vicenţiu D. Rădulescu[3,4,5]

**Abstract**
In this paper, we investigate the following fractional Sobolev critical Nonlinear Schrödinger coupled systems:

$$\begin{cases} (-\Delta)^s u = \mu_1 u + |u|^{2_s^*-2}u + \eta_1 |u|^{p-2}u + \gamma\alpha |u|^{\alpha-2}u|v|^{\beta} & \text{in } \mathbb{R}^N, \\ (-\Delta)^s v = \mu_2 v + |v|^{2_s^*-2}v + \eta_2 |v|^{q-2}v + \gamma\beta |u|^{\alpha}|v|^{\beta-2}v & \text{in } \mathbb{R}^N, \\ \|u\|_{L^2}^2 = m_1^2 \text{ and } \|v\|_{L^2}^2 = m_2^2, \end{cases}$$

where $(-\Delta)^s$ is the fractional Laplacian, $N > 2s$, $s \in (0, 1)$, $\mu_1, \mu_2 \in \mathbb{R}$ are unknown constants, which will appear as Lagrange multipliers, $2_s^*$ is the fractional Sobolev critical index, $\eta_1, \eta_2, \gamma, m_1, m_2 > 0$, $\alpha > 1$, $\beta > 1$, $p, q, \alpha + \beta \in (2 + 4s/N, 2_s^*]$. Firstly, if $p, q, \alpha + \beta < 2_s^*$, we obtain the existence of positive normalized solution when $\gamma$ is big enough. Secondly, if $p = q = \alpha + \beta = 2_s^*$, we show that nonexistence of positive normalized solution. The main ideas and methods of this paper are scaling transformation, classification discussion and concentration-compactness principle.

✉ Vicenţiu D. Rădulescu
radulescu@inf.ucv.ro

Jiabin Zuo
zuojiabin88@163.com

[1] College of Mathematics, Changchun Normal University, Changchun 130032, China

[2] School of Mathematics and Information Science, Guangzhou University, Guangzhou 510006, People's Republic of China

[3] Faculty of Applied Mathematics, AGH University of Science and Technology, 30-059 Kraków, Poland

[4] Simion Stoilow, Institute of Mathematics of the Romanian Academy, Calea Grivitei No. 21, 010702 Bucharest, Romania

[5] Department of Mathematics, University of Craiova, 200585 Craiova, Romania

## 1 Introduction and main result

The motivation for the problem studied in this article arises from finding stationary waves solutions of the following physical model:

$$\begin{cases} (-\Delta)^s \phi_1 = -i\dfrac{\partial \phi_1}{\partial t} + |\phi_1|^{2_s^*-2}\phi_1 + \eta_1|\phi_1|^{p-2}\phi_1 + \gamma\alpha|\phi_1|^{\alpha-2}\phi_1|\phi_2|^{\beta}, \\ (-\Delta)^s \phi_2 = -i\dfrac{\partial \phi_2}{\partial t} + |\phi_2|^{2_s^*-2}\phi_2 + \eta_2|\phi_2|^{q-2}\phi_2 + \gamma\beta|\phi_1|^{\alpha}|\phi_2|^{\beta-2}\phi_2, \\ \phi_j(x,t) \to 0 \text{ as } |x| \to \infty, \, j=1,2, \end{cases} \quad (1.1)$$

where $i$ represents the imaginary unit and $\phi_j = \phi_j(x,t) : \mathbb{R}^N \times \mathbb{R}^+ \to \mathbb{C}$ is the wave function of the $j$th ($j = 1, 2$) component, the mass of them represents the number of particles of each component in the mean-field models for binary mixtures of Bose-Einstein condensation, see [1–3] and references therein, $\eta_j$ and $\gamma$ denote the intraspecies and interspecies scattering lengths. The sign case for $\gamma$ determines whether the interaction between the states is attractive or repulsive, i.e. the interaction is attractive if $\gamma$ is positive, the interaction is repulsive if $\gamma$ is negative.

An important solution, known as *travelling* or *standing wave*, is characterized by *ansatz*

$$\phi_1(x,t) = e^{i\mu_1 t}u(x), \quad \phi_2(x,t) = e^{i\mu_2 t}v(x) \quad (1.2)$$

for two unknown functions $u, v : \mathbb{R}^N \to \mathbb{R}$, where $\mu_1, \mu_2 \in \mathbb{R}$. Because these solutions are very similar to each other and retain their mass over time, it makes sense to seek *prescribed $L^2$-norm solutions* (normalized solutions). Therefore, combining (1.1) and (1.2), we arrive at the following fractional single or double Sobolev critical Schrödinger system:

$$\begin{cases} (-\Delta)^s u = \mu_1 u + |u|^{2_s^*-2}u + \eta_1|u|^{p-2}u + \gamma\alpha|u|^{\alpha-2}u|v|^{\beta} \text{ in } \mathbb{R}^N, \\ (-\Delta)^s v = \mu_2 v + |v|^{2_s^*-2}v + \eta_2|v|^{q-2}v + \gamma\beta|u|^{\alpha}|v|^{\beta-2}v \text{ in } \mathbb{R}^N, \\ \|u\|_{L^2}^2 = m_1^2 \text{ and } \|v\|_{L^2}^2 = m_2^2, \end{cases} \quad (1.3)$$

where $N > 2s$, $s \in (0, 1)$, $\mu_1, \mu_2 \in \mathbb{R}$ are unknown constants, which will appear as Lagrange multipliers, $2_s^*$ is the fractional Sobolev critical index, $\eta_1, \eta_2, \gamma, m_1, m_2 > 0$, $\alpha > 1, \beta > 1, p, q, \alpha + \beta \in (2 + 4s/N, 2_s^*]$, and $(-\Delta)^s$ is the fractional Laplacican defined by

$$(-\Delta)^s u(x) = C(N,s) \lim_{\varepsilon \to 0^+} \int_{\mathbb{R}^N \setminus B_\varepsilon(x)} \frac{u(x) - u(y)}{|x-y|^{N+2s}} dxdy,$$

where

$$C(N,s) = \left(\int_{\mathbb{R}^N} \frac{1-\cos(\omega_1)}{\omega} d\omega\right)^{-1}.$$

For more information about this type of operator we refer to [4].

With regard to the double Sobolev critical Schrödinger coupled systems, many experts and scholars have conducted extensive and in-depth research on this whether it is integer order or fractional order with fixed $\mu_1$ and $\mu_2$. When $s \to 1$, Zou *et al.* [5] considered the existence and symmetry of positive ground states for a double critical coupled systems. Moreover, they studied the limit behavior of positive ground states for another kind of double critical Schrödinger system when the interaction is repulsive in [6]. In the case of $0 < s < 1$, Zou and Yin [7] proved the asymptotic behaviour and existence of the positive least energy solutions for k-coupled double critical systems driven by a fractional Laplace operator by means of the idea of induction. Yang [8] dealt with a class of fractional Laplacian doubly critical coupled systems, they gave sufficient conditions for the existence of weak solutions by establishing an embedding theorem. He et al. [9] investigated the existence of least energy solution with the help of the Nehari manifold.

However, as far as we know, few papers treat parameters $\mu_1$ and $\mu_2$ as Lagrange multipliers to study the normalized solution of double Sobolev critical problems.

In particular, when $s \to 1$, problem (1.3) becomes the form:

$$\begin{cases} -\Delta u = \mu_1 u + |u|^{2_s^*-2}u + \eta_1|u|^{p-2}u + \gamma\alpha|u|^{\alpha-2}u|v|^\beta & \text{in } \mathbb{R}^N, \\ -\Delta v = \mu_2 v + |v|^{2_s^*-2}v + \eta_2|v|^{q-2}v + \gamma\beta|u|^\alpha|v|^{\beta-2}v & \text{in } \mathbb{R}^N, \\ \|u\|_{L^2}^2 = m_1^2 \text{ and } \|v\|_{L^2}^2 = m_2^2. \end{cases} \quad (1.4)$$

Liu and Fang [10] studied the existence and nonexistence of positive normalized solution for equation (1.4) in the case of $p, q, \alpha + \beta < 2_s^*$ and $p = q = \alpha + \beta = 2_s^*$ respectively.

Furthermore, if let $N = 4$, $p = q$, $\alpha = \beta = 2$, the above problem reduces to the classical elliptic system:

$$\begin{cases} -\Delta u = \mu_1 u + u^3 + \eta_1|u|^{p-2}u + 2\gamma uv^2 & \text{in } \mathbb{R}^4, \\ -\Delta v = \mu_2 v + v^3 + \eta_2|v|^{p-2}v + 2\gamma\beta u^2 v & \text{in } \mathbb{R}^4, \\ \|u\|_{L^2}^2 = m_1^2 \text{ and } \|v\|_{L^2}^2 = m_2^2, \end{cases} \quad (1.5)$$

which was investigated by Zou et al. [11], they obtained the existence, nonexistence and asymptotic behavior of normalized ground state solutions for system (1.5) in different cases.

Of course, a great deal of work has focused on the normalized solution of integer-order nonlinear Schrödinger systems [12–18] or fractional-order single Schrödinger equations [19–26]. In particular, we highlight that Jeanjean and Lu [21] obtained the existence and multiplicity of normalized solutions and the asymptotic behavior of the ground state solution for a class of mass supercritical problems. The method developed in the present paper is inspired by the techniques introduced in [21]. We also point out that Soave [27] established existence and stability properties of ground state solutions for a Schrödinger equation with Sobolev critical exponent under three different assumptions, respectively mass subcritical, mass critical and mass supercritical.

So far, we have found only one paper [28] dealing with normalized solution of fractional Schrödinger coupled systems but only the subcritical case is considered. Hence it is natural to inquire what difficulties will appear if we consider single critical nonlinearity or even double critical nonlinearities.

Motivated by the work above, we will consider the existence for single critical fractional Schrödinger coupled systems and the nonexistence of normalized solutions for double critical fractional Schrödinger coupled systems. The main features of this paper is the existence of nonlocal operator and Sobolev critical nonlinearities, which makes dealing with compactness conditions more complicated. Our main methods and tools for solving these problems are scaling transformation, classification discussion and concentration-compactness principle.

A classical method for studying the normalized solution of system (1.3) is to look for critical points of the following $C^1$ functional

$$I(u,v) := \frac{1}{2}\left([u]_{H^s}^2 + [v]_{H^s}^2\right) - \frac{1}{2_s^*}\left(\|u\|_{2_s^*}^{2_s^*} + \|v\|_{2_s^*}^{2_s^*}\right) - \frac{\eta_1}{p}\|u\|_p^p$$
$$- \frac{\eta_2}{q}\|v\|_q^q - \gamma \int_{\mathbb{R}^N} |u|^\alpha |v|^\beta$$

constrained to the set

$$S_{m_1} \times S_{m_2} = \left\{(u,v) \in H^s(\mathbb{R}^N) \times H^s(\mathbb{R}^N) : \|u\|_2^2 = m_1^2, \ \|v\|_2^2 = m_2^2\right\},$$

where $H^s(\mathbb{R}^N)$ is the fractional Sobolev space defined by

$$H^s(\mathbb{R}^N) = \left\{u \in L^2(\mathbb{R}^N) \mid [u]_{H^s}^2 = \int_{\mathbb{R}^{2N}} \frac{u(x) - u(y)}{|x-y|^{\frac{N}{2}+s}} dxdy < \infty\right\},$$

whose norm is

$$\|u\| = (\|u\|_{L^2}^2 + [u]_{H^s}^2)^{\frac{1}{2}}.$$

For convenience, we use $\|u\|_p$ to represent the norm of Lebesgue space $L^p(\mathbb{R}^N)$ for $p \in [1, \infty)$. We write

$$\theta = \alpha + \beta, \ L^2\text{-critical exponent } \overline{p} = 2 + 4s/N,$$
$$H^s_{rad}(\mathbb{R}^N) = \{u \in H^s(\mathbb{R}^N) : u(x) = u(|x|)\}, \ S_m = \left\{u \in H^s(\mathbb{R}^N) : \|u\|_2^2 = m^2\right\},$$
$$S_{m_1,r} = S_{m_1} \cap H^s_{rad}(\mathbb{R}^N), \ S_{m_2,r} = S_{m_2} \cap H^s_{rad}(\mathbb{R}^N),$$
$$W_{m,r} = S_{m_1,r} \times S_{m_2,r}, \ W_r = H^s_{rad}(\mathbb{R}^N) \times H^s_{rad}(\mathbb{R}^N).$$

We now present our main results:

**Theorem 1.1** *Assume that $N > 2s$, $\alpha > 1$, $\beta > 1$, $p, q, \theta \in (\overline{p}, 2_s^*)$, and $\eta_1, \eta_2, \gamma, m_1, m_2 > 0$. Then there exists $\gamma^* = \gamma^*(m_1, m_2) > 0$ such that for any $\gamma \geq \gamma^*$, problem (1.3) admits a radial normalized solution $(\widetilde{u}, \widetilde{v})$. Moreover, $(\widetilde{u}, \widetilde{v})$ is a positive solution whose associated Lagrange multipliers $\mu_1$ and $\mu_2$ are negative.*

**Remark 1.1** The proof of Theorem 1.1 faces some difficulties and challenges. Firstly, strong convergence of sequences in $L^2(\mathbb{R}^N)$ space is difficult to obtain because the embeddings $H^s(\mathbb{R}^N) \hookrightarrow L^2(\mathbb{R}^N)$ and $H^s_{rad}(\mathbb{R}^N) = \{u \in H^s(\mathbb{R}^N) : u(x) = u(|x|)\} \hookrightarrow L^2(\mathbb{R}^N)$ are not compact. Secondly, the lack of compactness caused by Sobolev critical index makes verifying the Palais-Smale condition more complicated. Thirdly, the idea of classification discussion is going to be used since we don't infer which of the three indices $p, q, \theta$ is big and which is small.

**Remark 1.2** This result extends the partial result in [29] in some aspects. Compared with the local case of our result, this kind of system is studied by Mederski and Schino [29] under more generalized assumptions on the nonlinear terms.

**Theorem 1.2** *Assume that $N > 2s$, $\alpha > 1$, $\beta > 1$, $p = q = \theta = 2_s^*$, and $\eta_1, \eta_2, \gamma, m_1, m_2 > 0$. Then problem (1.3) has no positive normalized solution.*

**Remark 1.3** Theorems 1.1 and 1.2 seem to be the first results of the normalized solution for fractional Sobolev critical Schrödinger coupling systems.

The paper is organized as follows. Section 2 introduces relative results of scalar equations and some preliminaries, which play an important role in the proof of Palais-Smale condition. Section 3 proves Theorem 1.1 by using the methods of scaling transformation, classification discussion and concentration-compactness principle. Section 4 gives the proof of Theorem 1.2 with the help of Pohozaev identity.

## 2 Relevant results for scalar equations and preliminaries

In order to study the fractional critical Schrodinger coupling system, we first need to review related results of following scalar equations, i.e. $\gamma = 0$ in (1.3):

$$\begin{cases} (-\Delta)^s u = \mu_1 u + |u|^{2_s^*-2}u + \eta_1 |u|^{p-2}u & in \ \mathbb{R}^N, \\ u \in H^s(\mathbb{R}^N), \ \|u\|_2^2 = m_1^2, \end{cases} \quad (2.1)$$

which has been investigated in [25] by constraining on the Pohozaev manifold

$$\mathcal{P}_{m_1, \eta_1} = \{u \in S_{m_1} : [u]_{H^s}^2 - \|u\|_{2_s^*}^{2_s^*} - \eta_1 \xi_p \|u\|_p^p = 0\},$$

where

$$\xi_p = \frac{Np - 2N}{2ps}, \ \text{for any } p \in (2, 2_s^*]. \quad (2.2)$$

A standard way to get normalized solutions of (2.1) is to look for critical points for $C^1$ functional

$$I_{\eta_1}(u)|_{S_{m_1}} = \frac{1}{2}[u]^2_{H^s} - \frac{1}{2^*_s}\|u\|^{2^*_s}_{2^*_s} - \frac{\eta_1}{p}\|u\|^p_p,$$

As we all know, $\mathcal{P}_{m_1,\eta_1}$ contains every critical point of $I_{\eta_1}(u)|_{S_{m_1}}$, due to the Pohozaev identity (see [30, Proposition 4.1]).

It follows from [31] that there exists a best fractional critical Sobolev constant $\mathcal{S} > 0$ such that

$$\mathcal{S}\|u\|^2_{2^*_s} \le [u]^2_{H^s}, \text{ for all } u \in S_{m_1}, \tag{2.3}$$

which is famous fractional Sobolev inequality.

In order to prove our result in Sect. 3, we need to obtain the following monotonicity result of scalar equations, which is necessary in the proof of Lemma 3.7.

**Lemma 2.1** *Assume that $N > 2s$, $m_1, \eta_1 > 0$ and $p \in (\overline{p}, 2^*_s)$. Then, $m_1 \mapsto E^{\eta_1}_{m_1} = \inf_{\mathcal{P}_{m_1,\eta_1}} I_{\eta_1}(u) \in (0, \frac{s\mathcal{S}^{\frac{N}{2s}}}{N})$ is non-increasing in $(0, +\infty)$. where $\mathcal{S}$ is given in (2.3).*

**Proof** According to [25, Theorem 1.3], we obtain that $E^{\eta_1}_{m_1} \in (0, \frac{s\mathcal{S}^{\frac{N}{2s}}}{N})$, Similar to the proof of [19, Lemma 9] or [10, Lemma 2.1], we only need to make a small change to get non-increasing property of the function $E^{\eta_1}_{m_1}$ with respect to $m_1$, so we omit it. □

**Lemma 2.2** [22, Lemma 1.1] *There exists an optimal constant $C(N, p, s)$ such that for $p \in (2, 2^*_s)$,*

$$\|u\|^p_p \le C^p(N, p, s)[u]^{p\xi_p}_{H^s}\|u\|^{p(1-\xi_p)}_2, \quad \forall u \in H^s(\mathbb{R}^N), \tag{2.4}$$

*where $\xi_p$ is defined by (2.2).*

The above fractional Gagliardo-Nirenberg inequality plays an key role in the next series of proofs.

**Lemma 2.3** *Let $N > 2s$, $\alpha > 1$, $\beta > 1$, $\eta_1, \eta_2, \gamma > 0$, $p, q, \theta \in (2, 2^*_s)$, and $(u, v) \in W_r$ is a nonnegative solution of system (1.3). Then, it follows from $u \gneq 0$ that $\mu_1 < 0$; it follows from $v \gneq 0$ that $\mu_2 < 0$.*

**Proof** On account of $u \gneq 0$ and fulfills

$$(-\Delta)^s u = \mu_1 u + |u|^{2^*_s-2}u + \eta_1|u|^{p-2}u + \gamma\alpha|u|^{\alpha-2}u|v|^\beta \text{ in } \mathbb{R}^N,$$

thus we infer that $(-\Delta)^s u \ge 0$ if $\mu_1 \ge 0$. It follows from [32, Proposition 2.17] that $u \equiv 0$. This contradicts to the condition $u \gneq 0$, which means that $\mu_1 < 0$. Similarly, we also can get $\mu_2 < 0$ from $v \gneq 0$. □

**Lemma 2.4** *Let $N > 2s$, $\alpha > 1$, $\beta > 1$, $\theta \in [2, 2_s^*]$, $(u_n, v_n) \rightharpoonup (u, v)$ in $W_r$. Then, up to a subsequence*

$$\lim_{n \to \infty} \int_{\mathbb{R}^N} |u_n|^\alpha |v_n|^\beta - |u_n - u|^\alpha |v_n - v|^\beta - |u|^\alpha |v|^\beta = 0.$$

*Remark 2.1* Since the proof of Brézis-Lieb Lemma 2.4 is standard and classical, we would like to omit it, please refer to the literature [5, Lemma 2.3] for interested readers.

## 3 Proof of theorem 1.1

We will do a scaling transformation make the functional $I(\rho \star u, \rho \star v)$ satisfy the mountain pass geometry. For $(u, v) \in W$ and $\rho \in \mathbb{R}$, we let

$$(\rho \star u, \rho \star v) = \left( e^{\frac{N\rho}{2}} u(e^\rho x), e^{\frac{N\rho}{2}} v(e^\rho x) \right) \quad \text{for a.e. } x \in \mathbb{R}^N,$$

which comes from the inspiration of Jeanjean [33]. The results show that the original functional $I(u, v)$ and the transformed functional $\widetilde{I} = I(\rho \star u, \rho \star v)$ have the same mountain pass geometry and mountain pass level.

**Lemma 3.1** *Suppose that $(u, v) \in S_{m_1} \times S_{m_2}$ is arbitrary but fixed. Then we have the following conclusions:*

*(1) $[\rho \star u]^2_{H^s} + [\rho \star v]^2_{H^s} \to 0$ and $I(\rho \star u, \rho \star v) \to 0$ as $\rho \to -\infty$;*

*(2) $[\rho \star u]^2_{H^s} + [\rho \star v]^2_{H^s} \to +\infty$ and $I(\rho \star u, \rho \star v) \to -\infty$ as $\rho \to +\infty$.*

*Proof* Through simple calculations, we have

$$[\rho \star u]^2_{H^s} + [\rho \star v]^2_{H^s} = e^{2\rho s} \int \int_{\mathbb{R}^{2N}} \frac{|u(x) - u(y)|^2}{|x - y|^{N+2s}} dx dy$$
$$+ e^{2\rho s} \int \int_{\mathbb{R}^{2N}} \frac{|v(x) - v(y)|^2}{|x - y|^{N+2s}} dx dy$$
$$= e^{2\rho s} [u]^2_{H^s} + e^{2\rho s} [v]^2_{H^s}, \tag{3.1}$$

$$\|\rho \star u\|^\zeta_\zeta = e^{\frac{(\zeta - 2)N\rho}{2}} \|u\|^\zeta_\zeta,$$

$$\|\rho \star v\|^\zeta_\zeta = e^{\frac{(\zeta - 2)N\rho}{2}} \|v\|^\zeta_\zeta, \tag{3.2}$$

$$\gamma \int_{\mathbb{R}^N} |\rho \star u|^\alpha |\rho \star v|^\beta = \gamma e^{\frac{(\theta - 2)N\rho}{2}} \int_{\mathbb{R}^N} |u|^\alpha |v|^\beta. \tag{3.3}$$

From (3.1)–(3.3), $\zeta \geq 2$, and $\theta > 2$, we get

$$[\rho \star u]_{H^s}^2 \to 0, \quad [\rho \star v]_{H^s}^2 \to 0, \quad \text{as } \rho \to -\infty,$$
$$\|\rho \star u\|_\zeta^\zeta \to 0, \quad \|\rho \star v\|_\zeta^\zeta \to 0 \quad \text{as } \rho \to -\infty,$$
$$\gamma \int_{\mathbb{R}^N} |\rho \star u|^\alpha |\rho \star v|^\beta \to 0 \quad \text{as } \rho \to -\infty.$$

Thus, we have

$$I(\rho \star u, \rho \star v) = \frac{1}{2} \left([\rho \star u]_{H^s}^2 + [\rho \star v]_{H^s}^2\right) - \frac{1}{2_s^*} \left(\|\rho \star u\|_{2_s^*}^{2_s^*} + \|\rho \star v\|_{2_s^*}^{2_s^*}\right)$$
$$- \frac{\eta_1}{p} \|\rho \star u\|_p^p - \frac{\eta_2}{q} \|\rho \star v\|_q^q - \gamma \int_{\mathbb{R}^N} |\rho \star u|^\alpha |\rho \star v|^\beta \to 0 \quad \text{as } \rho \to -\infty,$$

which implies that (1) holds.

Again by (3.1), we obtain that $[\rho \star u]_{H^s}^2 + [\rho \star v]_{H^s}^2 \to +\infty$ as $\rho \to +\infty$. Furthermore,

$$I(\rho \star u, \rho \star v) = \frac{1}{2} \left([\rho \star u]_{H^s}^2 + [\rho \star v]_{H^s}^2\right) - \frac{1}{2_s^*} \left(\|\rho \star u\|_{2_s^*}^{2_s^*} + \|\rho \star v\|_{2_s^*}^{2_s^*}\right)$$
$$- \frac{\eta_1}{p} \|\rho \star u\|_p^p - \frac{\eta_2}{q} \|\rho \star v\|_q^q - \gamma \int_{\mathbb{R}^N} |\rho \star u|^\alpha |\rho \star v|^\beta \to -\infty \quad \text{as } \rho \to +\infty$$

since $p, q, \theta \in (2 + \frac{4s}{N}, 2_s^*)$, which also means that (2) holds. □

**Lemma 3.2** *There exists $A(m_1, m_2) > 0$ small enough such that*

$$0 < \sup_{u \in X} I(u, v) < \inf_{u \in Y} I(u, v),$$

*with*

$$X := \left\{(u, v) \in S_{m_1} \times S_{m_2} : [u]_{H^s}^2 + [v]_{H^s}^2 \leq A(m_1, m_2)\right\},$$
$$Y := \left\{(u, v) \in S_{m_1} \times S_{m_2} : [u]_{H^s}^2 + [v]_{H^s}^2 = 2A(m_1, m_2)\right\}.$$

**Proof** According to Lemma 2.2 and the Hölder inequality, for any $(u, v) \in S_{m_1} \times S_{m_2}$ we get that

$$\frac{\eta_1}{p}\|u\|_p^p \leq C^p(N,p,s,m_1,\eta_1)[u]_{H^s}^{p\xi_p} \leq C^p(N,p,s,m_1,\eta_1)\left([u]_{H^s}^2 + [v]_{H^s}^2\right)^{\frac{p\xi_p}{2}},$$
(3.4)

$$\frac{\eta_2}{q}\|u\|_q^q \leq C^q(N,q,s,m_2,\eta_2)[u]_{H^s}^{q\xi_q} \leq C^q(N,q,s,m_2,\eta_2)\left([u]_{H^s}^2 + [v]_{H^s}^2\right)^{\frac{q\xi_q}{2}},$$
(3.5)

$$\gamma \int_{\mathbb{R}^N} |u|^\alpha |v|^\beta \leq \gamma \|u\|_\theta^\alpha \|v\|_\theta^\beta \leq \gamma C(N,\alpha,\beta,s,m_1,m_2)\left([u]_{H^s}^2 + [v]_{H^s}^2\right)^{\frac{\theta\xi_\theta}{2}}.$$
(3.6)

On the one hand, if let $b = [u]_{H^s}^2 + [v]_{H^s}^2$ and $A > 0$ be arbitrary but fixed, then for any $(u,v) \in X$ such that $b \leq A$, it follows from (3.4)–(3.6) and (2.3) that

$$\begin{aligned} I(u,v) &\geq \frac{1}{2}b - \frac{1}{2_s^* \mathcal{S}^{\frac{2_s^*}{2}}} b^{\frac{2_s^*}{2}} - C^p(N,p,s,m_1,\eta_1) b^{\frac{p\xi_p}{2}} \\ &\quad - C^q(N,q,s,m_2,\eta_2) b^{\frac{q\xi_q}{2}} - \gamma C(N,\alpha,\beta,s,m_1,m_2) b^{\frac{\theta\xi_\theta}{2}} \\ &\geq \frac{1}{8}b > 0 \end{aligned}$$

for $A$ small enough, the scaling of inequality above takes advantage of this fact that $p\xi_p, q\xi_q, \theta\xi_\theta > 2$. On the other hand, for any $(u_1, v_1) \in Y$ and $(u_2, v_2) \in X$ such that $b_1 = [u_1]_{H^s}^2 + [v_1]_{H^s}^2 = 2A$ and $b_2 = [u_2]_{H^s}^2 + [v_2]_{H^s}^2 \leq A$, we obtain

$$\begin{aligned} I(u_1,v_1) - I(u_2,v_2) &\geq \frac{1}{2}(b_1 - b_2) - \frac{1}{2_s^* \mathcal{S}^{\frac{2_s^*}{2}}} b_1^{\frac{2_s^*}{2}} - C^p(N,p,s,m_1,\eta_1) b_1^{\frac{p\xi_p}{2}} \\ &\quad - C^q(N,q,s,m_2,\eta_2) b_1^{\frac{q\xi_q}{2}} - \gamma C(N,\alpha,\beta,s,m_1,m_2) b_1^{\frac{\theta\xi_\theta}{2}} \\ &\geq \frac{1}{2}A - \frac{1}{2_s^* \mathcal{S}^{\frac{2_s^*}{2}}}(2A)^{\frac{2_s^*}{2}} - C^p(N,p,s,m_1,\eta_1)(2A)^{\frac{p\xi_p}{2}} \\ &\quad - C^q(N,q,s,m_2,\eta_2)(2A)^{\frac{q\xi_q}{2}} - \gamma C(N,\alpha,\beta,s,m_1,m_2)(2A)^{\frac{\theta\xi_\theta}{2}} \\ &\geq \frac{1}{8}b \end{aligned}$$

for $A$ sufficiently small. So, we can pick $A$ small enough for the inequality in Lemma 3.2 to be true. □

Now we have obtained that the geometry of mountain pass, then we give the minimax picture: Because Zhang [25] had proved the $\bar{u} := u_{m_1,\eta_1}$ is a ground state of (2.1) involving parameters $p, \eta_1, m_1$ and also $\bar{v} := u_{m_2,\eta_2}$ is a ground state of (2.1) involving parameters $p, \eta_2, m_2$. Therefore, we fix $(\bar{u}, \bar{v}) \in W_{m,r}$, according to Lemma

3.1 and Lemma 3.2, there exist two numbers $\rho_1 \ll -1 < 0 < 1 \ll \rho_2$ such that

$$e^{2\rho_1 s}[\overline{u}]_{H^s}^2 + e^{2\rho_1 s}[\overline{v}]_{H^s}^2 < \frac{A(m_1, m_2)}{2}, \quad I(\rho_1 \star \overline{u}, \rho_1 \star \overline{v}) > 0,$$
$$e^{2\rho_2 s}[\overline{u}]_{H^s}^2 + e^{2\rho_2 s}[\overline{v}]_{H^s}^2 > 2A(m_1, m_2), \quad I(\rho_2 \star \overline{u}, \rho_2 \star \overline{v}) \leq 0.$$

We define the path

$$\Gamma = \{\chi \in C([0, 1], W_{m,r}) : \chi(0) = (\rho_1 \star \overline{u}, \rho_1 \star \overline{v}), \chi(1) = (\rho_2 \star \overline{u}, \rho_2 \star \overline{v})\},$$

then $\Gamma$ is not empty. In fact, let $\chi_0(t) := ([(1-t)\rho_1 + t\rho_2] \star \overline{u}, [(1-t)\rho_1 + t\rho_2] \star \overline{v})$, obviously, $\chi_0(t) \in W_{m,r}$ and $\chi_0(t) \in \Gamma$.

Letting

$$c_\gamma(m_1, m_2) := \inf_{\chi \in \Gamma} \max_{t \in (0,1]} I(\chi(t)),$$

clearly, $c_\gamma(m_1, m_2) > 0$, then we have the following asymptotic behavior of critical values:

**Lemma 3.3** $\lim_{\gamma \to +\infty} c_\gamma(m_1, m_2) = 0.$

**Proof** Fix $(u_0, v_0) \in W_{m,r}$ and it follows from the path $\chi_0(t) := ([(1-t)\rho_1 + t\rho_2] \star u_0, [(1-t)\rho_1 + t\rho_2] \star v_0)$ that

$$c_\gamma(m_1, m_2) \leq \max_{t \in [0,1]} I(\chi_0(t))$$
$$\leq \max_{l \geq 0} \left\{ \frac{1}{2} l^2 \left([u_0]_{H^s}^2 + [v_0]_{H^s}^2\right) - \gamma l^{\frac{N\theta - 2N}{2s}} \int_{\mathbb{R}^N} |u_0|^\alpha |v_0|^\beta \right\}.$$

Let $C_1 = [u_0]_{H^s}^2 + [v_0]_{H^s}^2$ and $C_2 = \int_{\mathbb{R}^N} |u_0|^\alpha |v_0|^\beta$, we discuss the maximum value of the following function

$$g(l) = \frac{1}{2} C_1 l^2 - \gamma C_2 l^{\frac{N\theta - 2N}{2s}}, \text{ for any } l \geq 0.$$

Letting

$$g'(l) = C_1 l - \left(\frac{N\theta - 2N}{2s}\right) \gamma C_2 l^{\frac{N\theta - 2N - 2s}{2s}} = 0,$$

we can obtain the maximum point of $g(l)$, that is

$$l_{max} = \left(\frac{2s C_1}{(N\theta - 2N)\gamma C_2}\right)^{\frac{2s}{N\theta - 2N - 4s}}.$$

Thus,

$$\max_{l \geq 0} \left\{ \frac{1}{2}l^2 \left([u_0]_{H^s}^2 + [v_0]_{H^s}^2\right) - \gamma l^{\frac{N\theta-2N}{2s}} \int_{\mathbb{R}^N} |u_0|^\alpha |v_0|^\beta \right\}$$

$$= \frac{1}{2} \left(\frac{2sC_1}{(N\theta - 2N)\gamma C_2}\right)^{\frac{4s}{N\theta-2N-4s}} C_1 - \gamma \left(\frac{2sC_1}{(N\theta - 2N)\gamma C_2}\right)^{\frac{N\theta-2N}{N\theta-2N-4s}} C_2$$

$$\leq \frac{1}{2} \left(\frac{2sC_1}{(N\theta - 2N)\gamma C_2}\right)^{\frac{4s}{N\theta-2N-4s}} C_1.$$

As a result, there exists $C > 0$ that don't depend on $\gamma > 0$ such that

$$c_\gamma(m_1, m_2) \leq C \left(\frac{1}{\gamma}\right)^{\frac{4s}{N\theta-2N-4s}} \to 0 \text{ as } \gamma \to \infty$$

thanks to $\theta > 2 + 4s/N$. We complete the proof of Lemma 3.3 now. □

To further prove our main results, we need to define Pohozaev manifold $\mathcal{P}_{m_1,m_2,\mu_1,\mu_2}$ of vector equations (1.3) in the same way as scalar equations (2.1).

$$\mathcal{P}_{m_1,m_2,\mu_1,\mu_2} = \left\{(u, v) \in S_{m_1} \times S_{m_2} : P(u, v) = 0\right\},$$

where

$$P(u, v) := [u]_{H^s}^2 + [v]_{H^s}^2 - \|u\|_{2_s^*}^{2_s^*} - \|v\|_{2_s^*}^{2_s^*} - \eta_1 \xi_p \|u\|_p^p$$
$$- \eta_2 \xi_q \|v\|_q^q - \gamma \theta \xi_\theta \int_{\mathbb{R}^N} |u|^\alpha |v|^\beta.$$

Notice that $I(|u|, |v|) = I(u, v)$ and $P(|u|, |v|) = P(u, v)$ for any $(u, v) \in W_{m,r}$. Similar to the argument in [33, Proposition 2.2] with minor changes, replace $(0, \chi_n) \in \overline{\Gamma}$ with $(0, |\chi_n|) \in \overline{\Gamma}$ if necessary, where

$$\overline{\Gamma} := \left\{\overline{\chi} \in ([0, 1], W_{m,r}) : \overline{\chi}(0) = (0, (\rho_1 \star \overline{u}, \rho_1 \star \overline{v})), \overline{\chi}(1) = (0, (\rho_2 \star \overline{u}, \rho_2 \star \overline{v}))\right\}.$$

Then $\{(u_n, v_n)\} \subset W_{m,r}$ is a Palais-Smale sequence for $I(u, v)$ at level $c_\gamma(m_1, m_2)$, i.e.

$$I(u_n, v_n)) \to c_\gamma(m_1, m_2) \text{ as } n \to \infty \text{ and } I'(u_n, v_n) \to 0 \text{ as } n \to \infty.$$

Moreover,

$$u_n^- \to 0, \ v_n^- \to 0 \text{ a.e. in } \mathbb{R}^N \text{ as } n \to \infty.$$

**Lemma 3.4** *Suppose that $\{(u_n, v_n)\} \subset W_{m,r}$ is a Palais-Smale sequence for $I(u, v)$. Then $\lim_{n\to\infty} P(u_n, v_n) = 0$.*

**Proof** Let

$$(\rho_n \star u_n, \rho_n \star v_n) = (x_n, y_n).$$

By direct calculations, we claim that

$$\frac{d}{d\rho_n} I((-\rho_n)\star x_n, (-\rho_n)\star y_n) = -sP((-\rho_n)\star x_n, (-\rho_n)\star y_n) = -sP(u_n, v_n).$$

Since $I'(u_n, v_n) \to 0$ as $n \to \infty$, again from [22, Proposition 5.4(2)], we get that $\lim_{n\to\infty} P(u_n, v_n) = 0$. □

**Lemma 3.5** *If $\{(u_n, v_n)\} \subset W_{m,r}$ is a Palais-Smale sequence for $I(u, v)$, then $\{(u_n, v_n)\}$ is bounded in $W_{m,r}$.*

**Proof** Observing that $\xi_p p, \xi_q q, \xi_\theta \theta > 2$ due to the fact that $p, q, \theta > 2 + 4s/N$. It follows from Lemma 3.4 that

$$I(u_n, v_n) = \frac{\eta_1}{2p}(\xi_p p - 2)\|u_n\|_p^p + \frac{\eta_2}{2q}(\xi_q q - 2)\|v_n\|_q^q + \frac{\gamma}{2}(\xi_\theta \theta - 2)\int_{\mathbb{R}^N} |u|^\alpha |v|^\beta$$
$$+ \frac{s}{N}\|u_n\|_{2_s^*}^{2_s^*} + \frac{s}{N}\|v_n\|_{2_s^*}^{2_s^*} + o(1).$$

Since $I(u_n, v_n)$ is bounded, thus we conclude that sequences $\{\|u_n\|_p^p\}$, $\{\|v_n\|_q^q\}$, $\{\int_{\mathbb{R}^N} |u_n|^\alpha |v_n|^\beta\}$, $\{\|u_n\|_{2_s^*}^{2_s^*}\}$, and $\{\|v_n\|_{2_s^*}^{2_s^*}\}$ are all bounded. Again by Lemma 3.4, we infer that $\{[u_n]_{H^s}^2\}$ and $\{[v_n]_{H^s}^2\}$ are also bounded, as claimed. □

In view of Lemma 3.5, there exists a nonnegative $(\widetilde{u}, \widetilde{v}) \in W_r$ such that, up to a subsequence,

$$\begin{aligned}
(u_n, v_n) &\rightharpoonup (\widetilde{u}, \widetilde{v}) \text{ in } W_r; \\
(u_n, v_n) &\rightharpoonup (\widetilde{u}, \widetilde{v}) \text{ in } L^{2_s^*}(\mathbb{R}^N) \times L^{2_s^*}(\mathbb{R}^N); \\
(u_n, v_n) &\to (\widetilde{u}, \widetilde{v}) \text{ in } L^p(\mathbb{R}^N) \times L^q(\mathbb{R}^N); \\
(u_n, v_n) &\to (\widetilde{u}, \widetilde{v}) \text{ in } L^\theta(\mathbb{R}^N) \times L^\theta(\mathbb{R}^N); \\
(u_n, v_n) &\to (\widetilde{u}, \widetilde{v}) \text{ a.e. in } \mathbb{R}^N
\end{aligned} \quad (3.7)$$

as $n \to \infty$. Since $\{(u_n, v_n)\} \subset S_{m_1} \times S_{m_2}$ is a Palais-Smale sequence for $I(u, v)$, on basis of the Lagrange multipliers rule, there exists a sequence $\{(\mu_1^n, \mu_2^n)\} \subset \mathbb{R} \times \mathbb{R}$ such that

$$I'(u_n, v_n) + \mu_1^n(u_n, 0) + \mu_2^n(0, v_n) \to 0 \text{ in } W_r \text{ as } n \to \infty. \quad (3.8)$$

Next, we take $(u_n, 0)$ and $(0, v_n)$ as test functions in (3.7), it follows from the proof of [25, Proposition 2.2] that $\{(\mu_1^n, \mu_2^n)\}$ is bounded in $\mathbb{R}$. Therefore up to a subsequence $(\mu_1^n, \mu_2^n) \to (\mu_1, \mu_1) \in \mathbb{R} \times \mathbb{R}$.

**Lemma 3.6** *There exists $\gamma^* = \gamma^*(m_1, m_2) > 0$ big enough, such that for any $\gamma \geq \gamma^*$, $(u_n, v_n) \to (\widetilde{u}, \widetilde{v})$ in $L^{2^*_s}(\mathbb{R}^N) \times L^{2^*_s}(\mathbb{R}^N)$ and $\widetilde{u}, \widetilde{v} \neq 0$.*

*Proof* We first claim that $u_n \to \widetilde{u}$ in $L^{2^*_s}(\mathbb{R}^N)$. In fact, according to the concentration-compactness principle in [34], we know that there exist two nonnegative measures $\omega$, $\nu$ and a (at most countable) index set $J$ such that $|(-\Delta)^{\frac{s}{2}} u_n|^2 \rightharpoonup \omega$ in $\mathcal{M}(\mathbb{R}^N)$, $|u_n|^{2^*_s} \rightharpoonup \nu$ in $\mathcal{M}(\mathbb{R}^N)$ and

$$\begin{cases} \omega \geq |(-\Delta)^{\frac{s}{2}} \widetilde{u}|^2 + \sum_{j \in J} \omega_j \delta_{x_j}, & \omega_j \geq 0, \\ \nu = |\widetilde{u}|^{2^*_s} + \sum_{j \in J} \nu_j \delta_{x_j}, & \nu_j \geq 0, \\ \nu_j \leq \mathcal{S}^{-\frac{2^*_s}{2}} \omega_j^{\frac{2^*_s}{2}}, & \forall j \in J, \end{cases} \quad (3.9)$$

where $x_j$ is the different point in $\mathbb{R}^N$, $\delta_{x_j}$ denotes the Dirac measure at $x_j$, $\mathcal{S}$ is best Sobolev constant in (2.3). Furthermore, it is possible to lose mass at infinity. i.e.,

$$\begin{cases} \limsup_{n \to \infty} \int_{\mathbb{R}^N} |(-\Delta)^{\frac{s}{2}} u_n|^2 dx = \int_{\mathbb{R}^N} d\omega + \omega_\infty, \\ \limsup_{n \to \infty} \int_{\mathbb{R}^N} |u_n|^{2^*_s} dx = \int_{\mathbb{R}^N} d\nu + \nu_\infty, \\ \nu_\infty \leq \mathcal{S}^{-\frac{2^*_s}{2}} \omega_\infty^{\frac{2^*_s}{2}}, \end{cases} \quad (3.10)$$

where

$$\omega_\infty = \lim_{R \to \infty} \limsup_{n \to \infty} \int_{|x| > R} |(-\Delta)^{\frac{s}{2}} u_n|^2 dx,$$
$$\nu_\infty = \lim_{R \to \infty} \limsup_{n \to \infty} \int_{|x| > R} |u_n|^{2^*_s} dx.$$

*Case 1* Take $\psi_\epsilon(x) \in C_0^\infty(\mathbb{R}^N)$ be a cut-off function such that

$$\begin{cases} \psi_\epsilon(x) \equiv 1 \text{ in } B_\epsilon(x_j), \\ \psi_\epsilon(x) \equiv 0 \text{ in } B_{2\epsilon}^c(x_j), \\ 0 \leq \psi_\epsilon(x) \leq 1, \end{cases} \quad (3.11)$$

where $B_\epsilon(x_j)$ represents the small ball with radius $\epsilon$ and center $x_j$. By Lemma 3.5, we note that $\{\psi_\epsilon(x) u_n\}$ is bounded in $H^s(\mathbb{R}^N)$. Next, again take $\{(\psi_\epsilon(x) u_n, 0)\}$ as a test function in (3.8) and letting $\epsilon \to 0$, we obtain that

$$\lim_{\epsilon \to 0} \lim_{n \to \infty} \langle I'(u_n, v_n) - \mu_1^n (u_n, 0) - \mu_2^n (0, v_n), (\psi_\epsilon(x) u_n, 0) \rangle = 0. \quad (3.12)$$

From (3.7), the Hölder inequality and the absolute continuity of the Lebesgue integral, we have

$$\lim_{\epsilon \to 0} \lim_{n \to \infty} \int_{\mathbb{R}^N} \mu_1^n u_n^2 \psi_\epsilon dx = 0,$$

$$\lim_{\epsilon \to 0} \lim_{n \to \infty} \int_{\mathbb{R}^N} |u_n|^p \psi_\epsilon dx = \lim_{\epsilon \to 0} \int_{\mathbb{R}^N} |\widetilde{u}|^p \psi_\epsilon dx = 0, \quad (3.13)$$

$$\lim_{\epsilon \to 0} \lim_{n \to \infty} \int_{\mathbb{R}^N} |u_n|^\alpha |v_n|^\beta \psi_\epsilon dx = \lim_{\epsilon \to 0} \int_{\mathbb{R}^N} |\widetilde{u}|^\alpha |\widetilde{v}|^\beta \psi_\epsilon dx = 0.$$

According to (3.12) and (3.13), we deduce that

$$\lim_{\epsilon \to 0} \lim_{n \to \infty} \int_{\mathbb{R}^N} |(-\Delta)^{\frac{s}{2}} u_n|^2 \psi_\epsilon dx = \lim_{\epsilon \to 0} \lim_{n \to \infty} \int_{\mathbb{R}^N} |u_n|^{2_s^*} \psi_\epsilon dx,$$

which means that

$$\lim_{\epsilon \to 0} \int_{\mathbb{R}^N} \psi_\epsilon d\omega = \lim_{\epsilon \to 0} \int_{\mathbb{R}^N} \psi_\epsilon d\nu. \quad (3.14)$$

Therefore it follows from (3.9) and (3.14) that $\nu_j \geq \omega_j$, thereby

$$\text{either } \omega_j = 0 \text{ or } \omega_j \geq \mathcal{S}^{\frac{N}{2s}} \text{ for } j \in J, \quad (3.15)$$

which means that $J$ is a finite set.

*Case 2* Take $\varphi \in C_0^\infty(\mathbb{R}^N)$ be another cut-off function with

$$\begin{cases} 0 \leq \varphi \leq 1, \\ \varphi \equiv 0 \text{ in } B_{\frac{1}{2}}(0), \\ \varphi \equiv 1 \text{ in } \mathbb{R}^N \setminus B_1(0). \end{cases} \quad (3.16)$$

For any $R$, let $\varphi_R(x) = \varphi(\frac{x}{R})$, by Lemma 3.5, we also have that $\{\varphi_R(x) u_n\}$ is bounded in $H^s(\mathbb{R}^N)$. Similarly, take $\{(\varphi_R(x) u_n, 0)\}$ as a test function in (3.8) and letting $R \to \infty$, we also obtain that

$$\lim_{R \to \infty} \lim_{n \to \infty} \langle I'(u_n, v_n) - \mu_1^n(u_n, 0) - \mu_2^n(0, v_n), (\varphi_R(x) u_n, 0) \rangle = 0. \quad (3.17)$$

Similarly to the proof of [35, Lemma 3.3], we can also get

$$\lim_{R \to \infty} \lim_{n \to \infty} \int_{\mathbb{R}^N} \mu_1^n u_n^2 \varphi_R dx = 0,$$

$$\lim_{R \to \infty} \lim_{n \to \infty} \int_{\mathbb{R}^N} |u_n|^p \varphi_R dx = \lim_{R \to \infty} \int_{\mathbb{R}^N} |\widetilde{u}|^p \varphi_R dx = 0, \quad (3.18)$$

$$\lim_{R \to \infty} \lim_{n \to \infty} \int_{\mathbb{R}^N} |u_n|^\alpha |v_n|^\beta \varphi_R dx = \lim_{R \to \infty} \int_{\mathbb{R}^N} |\widetilde{u}|^\alpha |\widetilde{v}|^\beta \varphi_R dx = 0.$$

By (3.17) and (3.18), we also have

$$\lim_{R\to\infty}\lim_{n\to\infty}\int_{\mathbb{R}^N}|(-\Delta)^{\frac{s}{2}}u_n|^2\varphi_R dx = \lim_{R\to\infty}\lim_{n\to\infty}\int_{\mathbb{R}^N}|u_n|^{2^*_s}\varphi_R dx,$$

that is

$$\omega_\infty = \nu_\infty. \qquad (3.19)$$

Again by (3.10), we have

$$\text{either } \omega_\infty = 0 \text{ or } \omega_\infty \geq \mathcal{S}^{\frac{N}{2s}}. \qquad (3.20)$$

For the rest of the proof, we will only prove Case 1, because Case 2 and Case 1 are almost exactly the same.

If for any $j \in J$, $\omega_j = 0$, then we have that $\nu_j = 0$ since (3.9), thereby $|u_n|^{2^*_s} \to |\tilde{u}|^{2^*_s}$. By the Brézis-Lieb lemma [36], we conclude that $u_n \to \tilde{u}$ in $L^{2^*_s}(\mathbb{R}^N$, as claimed.

If instead $\omega_j \geq \mathcal{S}^{\frac{N}{2s}}$ for some $j \in J$. In view of $I(u_n, v_n) \to c_\gamma(m_1, m_2)$, Lemma 3.4 and (3.11), we adopt the method of categorical discussion:

(1) If $\theta = \min\{p, q, \theta\}$. It follows from Lemma 3.3 that there exists a positive constant $\gamma^1$ big enough, such that $c_\gamma(m_1, m_2) < \left(\frac{1}{2} - \frac{1}{\theta\xi_\theta}\mathcal{S}^{\frac{N}{2s}}\right)$ for any $\gamma \geq \gamma'$. By (3.9), we have

$$\left(\frac{1}{2} - \frac{1}{\theta\xi_\theta}\right)\mathcal{S}^{\frac{N}{2s}} > c_\gamma(m_1, m_2) = \lim_{n\to\infty}\left(I(u_n, v_n) - \frac{1}{\theta\xi_\theta}P(u_n, v_n)\right)$$
$$\geq \left(\frac{1}{2} - \frac{1}{\theta\xi_\theta}\right)\limsup_{n\to\infty}\int_{\mathbb{R}^N}|(-\Delta)^{\frac{s}{2}}u_n|^2\psi_\epsilon dx$$
$$= \left(\frac{1}{2} - \frac{1}{\theta\xi_\theta}\right)\int_{\mathbb{R}^N}\psi_\epsilon d\omega$$
$$\geq \left(\frac{1}{2} - \frac{1}{\theta\xi_\theta}\right)\omega_j \geq \left(\frac{1}{2} - \frac{1}{\theta\xi_\theta}\right)\mathcal{S}^{\frac{N}{2s}},$$

which is a contradiction.

(2) If $p = \min\{p, q, \theta\}$. Similar to (**1**), there is also a constant $\gamma''$ big enough, such that $c_\gamma(m_1, m_2) < \left(\frac{1}{2} - \frac{1}{p\xi_p}\mathcal{S}^{\frac{N}{2s}}\right)$ for any $\gamma \geq \gamma''$. We also obtain

$$c_\gamma(m_1, m_2) = \lim_{n\to\infty}\left(I(u_n, v_n) - \frac{1}{p\xi_p}P(u_n, v_n)\right)$$
$$\geq \left(\frac{1}{2} - \frac{1}{p\xi_p}\right)\omega_j \geq \left(\frac{1}{2} - \frac{1}{p\xi_p}\right)\mathcal{S}^{\frac{N}{2s}},$$

which contradicts our hypothesis.

(3) If $q = \min\{p, q, \theta\}$. Analogously as (1) and (2), we can also get our conclusion, which we omit here.

To sum up, there exists a bigger positive constant $\gamma^*(m_1, m_2)$ such that $\omega_j = 0 = \nu_j$ for any $j \in J$ and $\gamma \geq \gamma^*$. Therefore, for $\gamma \geq \gamma^*$, we have $u_n \to \widetilde{u}$ in $L^{2_s^*}(\mathbb{R}^N)$. The proof for $v_n \to \widetilde{v}$ in $L^{2_s^*}(\mathbb{R}^N)$ is similar.

Finally, we claim that $\widetilde{u}, \widetilde{v} \neq 0$. In fact, if not, we have $(\widetilde{u}, \widetilde{v}) = (0, 0)$. According to (3.6), (3.7), $u_n \to \widetilde{u}, v_n \to \widetilde{v}$ in $L^{2_s^*}(\mathbb{R}^N)$ and Lemma 3.4, we have

$$\lim_{n \to \infty} ([u_n]_{H^s}^2 + [v_n]_{H^s}^2) = 0. \tag{3.21}$$

It follows from (3.21) that $c_\gamma(m_1, m_2) = \lim_{n \to \infty} I(u_n, v_n) = 0$, which is impossible since $c_\gamma(m_1, m_2) > 0$. As a result, we end the proof. □

Set $c_0(m_1, 0) := E_{m_1}^{\eta_1}$ and $c_0(0, m_2) := E_{m_2}^{\eta_2}$.

**Remark 3.1** For any $m_1, m_2 > 0$, according to Lemma 3.3, if necessary, to choose a bigger $\gamma^*$ such that $c_\gamma(m_1, m_2) < \min\{c_0(m_1, 0), c_0(0, m_2)\}$ for any $\gamma \geq \gamma^*$.

**Lemma 3.7** *Assume that $c_\gamma(m_1, m_2) < \min\{c_0(m_1, 0), c_0(0, m_2)\}$, then $(u_n, v_n) \to (\widetilde{u}, \widetilde{v})$ in $W_r$. Moreover, $(\widetilde{u}, \widetilde{v}) \in W_r$ is a positive normalized solution for system (1.3) associated with $\mu_1, \mu_2 < 0$.*

**Proof** From Lemma 3.6, we know that $\widetilde{u}, \widetilde{v} \neq 0$. Next, we are going to classify it into two cases and prove by contradiction:

*Case 1* $\widetilde{u} \gneq 0, \widetilde{v} \equiv 0$. From the strong maximum principle for fractional Laplace operators [32, Proposition 2.17], we have that $\widetilde{u} > 0$. It follows from Lemma 2.1 and [25, Theorem 1.3] that $\widetilde{u}$ is a positive radial symmetric solution for problem (2.1) with parameters $p, \eta_1, m_1'$ where $m_1' = \|\widetilde{u}\| \leq m_1^2$, and $c_0(m_1, 0) \leq c_0(m_1', 0) \leq I(\widetilde{u}, 0)$.

In view of (3.7), Lemmas 2.4, 3.6, and two kinds of Brézis-Lieb Lemmas [36–38], we obtain that

$$0 = \lim_{n \to \infty} P(u_n, v_n) = \lim_{n \to \infty} P(u_n - \widetilde{u}, v_n) + P(\widetilde{u}, 0)$$
$$= \lim_{n \to \infty} [u_n - \widetilde{u}]_{H^s}^2 + [v_n]_{H^s}^2,$$

and

$$c_\gamma(m_1, m_2) = \lim_{n \to \infty} I(u_n, v_n) = \lim_{n \to \infty} I(u_n - \widetilde{u}, v_n) + I(\widetilde{u}, 0)$$
$$\geq \frac{1}{2} \lim_{n \to \infty} [u_n - \widetilde{u}]_{H^s}^2 + [v_n]_{H^s}^2 + c_0(m_1, 0)$$
$$= c_0(m_1, 0).$$

This contradicts our hypothesis.

*Case 2* $\widetilde{u} \equiv 0, \widetilde{v} \gneq 0$. Similarly to **Case 1**, $\widetilde{v}$ is a positive radial symmetric solution for problem (2.1) with parameters $q, \eta_2, m_2'$ where $m_2' = \|\widetilde{v}\| \leq m_2^2$, and $c_0(0, m_2) \leq$

$c_0(0, m_2') \leq I(0, \widetilde{v})$. We also have that

$$c_\gamma(m_1, m_2) = \lim_{n \to \infty} I(u_n, v_n) = \lim_{n \to \infty} I(u_n, v_n - \widetilde{v}) + I(0, \widetilde{v})$$
$$\geq c_0(m_1, 0),$$

which is still a contradiction. So, we obtain that $\widetilde{u} \not\equiv 0$ and $\widetilde{v} \not\equiv 0$. Therefore, $\mu_1, \mu_2 < 0$ follows from Lemma 2.3 and (3.7) and (3.8).

On basis of (3.7) and (3.8), Lemma 3.4 and the boundedness of sequence $\{(\mu_1^n, \mu_2^n)\}$, we have

$$\mu_1 m_1^2 + \mu_2 m_2^2 = \lim_{n \to \infty} (\mu_1^n \|u_n\|_2^2 + \mu_2^n \|v_n\|_2^2)$$
$$= \lim_{n \to \infty} \left( (\xi_p - 1)\eta_1 \|u_n\|_p^p + (\xi_q - 1)\eta_2 \|v_n\|_q^q + (\xi_\theta - 1)\gamma\theta \int_{\mathbb{R}^N} |u_n|^\alpha |v_n|^\beta \right)$$
$$= (\xi_p - 1)\eta_1 \|\widetilde{u}\|_p^p + (\xi_q - 1)\eta_2 \|\widetilde{v}\|_q^q + (\xi_\theta - 1)\gamma\theta \int_{\mathbb{R}^N} |\widetilde{u}|^\alpha |\widetilde{v}|^\beta$$
$$= \mu_1 \|\widetilde{u}\|_2^2 + \mu_2 \|\widetilde{v}\|_2^2,$$

which means that

$$\mu_1(\|\widetilde{u}\|_2^2 - m_1^2) + \mu_2(\|\widetilde{v}\|_2^2 - m_2^2) = 0.$$

Thus,

$$\|\widetilde{u}\|_2^2 = m_1^2, \quad \|\widetilde{v}\|_2^2 = m_2^2, \text{ i.e., } (u_n, v_n) \to (\widetilde{u}, \widetilde{v}) \text{ in } L^2 \times L^2.$$

The proof is now complete. □

Again by (3.7) and (3.8) and Lemma 3.4, we have

$$\lim_{n \to \infty} \left( [u_n]_{H^s}^2 + \mu_1 \|u_n\|_2^2 \right) = \lim_{n \to \infty} \left( \|u_n\|_{2_s^*}^{2_s^*} + \eta_1 \|u_n\|_p^p + \gamma\alpha \int_{\mathbb{R}^N} |u_n|^\alpha |v_n|^\beta \right)$$
$$= \|\widetilde{u}\|_{2_s^*}^{2_s^*} + \eta_1 \|\widetilde{u}\|_p^p + \gamma\alpha \int_{\mathbb{R}^N} |\widetilde{u}|^\alpha |\widetilde{v}|^\beta$$
$$= [\widetilde{u}]_{H^s}^2 + \mu_1 \|\widetilde{u}\|_2^2,$$

which implies that $\|u_n\|_{W_r} \to \|\widetilde{u}\|_{W_r}$ as $n \to \infty$. Similar to the above argument, we also obtain

$$\lim_{n \to \infty} \left( [v_n]_{H^s}^2 + \mu_1 \|v_n\|_2^2 \right) = [\widetilde{v}]_{H^s}^2 + \mu_1 \|\widetilde{v}\|_2^2,$$

which also means that $\|v_n\|_{W_r} \to \|\widetilde{v}\|_{W_r}$ as $n \to \infty$. Therefore, we conclude that $(u_n, v_n) \to (\widetilde{u}, \widetilde{v})$ in $W_r \times W_r$. The proof of Lemma 3.7 is complete.

**Proof of Theorem 1.1** It follows from Lemmas 3.5–3.7 that the proof of Theorem 1.1 is complete. □

## 4 Proof of theorem 1.2

**Lemma 4.1** *Assume that $N > 2s$, $\alpha > 1$, $\beta > 1$, $p = q = \alpha + \beta + \theta = 2_s^*$, and $\eta_1, \eta_2, \gamma, m_1, m_2 > 0$. Then the following system*

$$\begin{cases} (-\Delta)^s u = \mu_1 u + |u|^{2_s^*-2}u + \eta_1|u|^{2_s^*-2}u + \gamma\alpha|u|^{\alpha-2}u|v|^\beta & \text{in } \mathbb{R}^N, \\ (-\Delta)^s v = \mu_2 v + |v|^{2_s^*-2}v + \eta_2|v|^{2_s^*-2}v + \gamma\beta|u|^\alpha|v|^{\beta-2}v & \text{in } \mathbb{R}^N, \\ \|u\|_{L^2}^2 = m_1^2 \text{ and } \|v\|_{L^2}^2 = m_2^2, \quad u, v \in H^s(\mathbb{R}^N) \end{cases} \quad (4.1)$$

*has no positive normalized solution.*

**Proof** We will use the proof by contradiction. If $(u, v)$ is a positive solution for the system (4.1) with some $\mu_1, \mu_2 \in \mathbb{R}$. According to Lemma 2.3, we know that $\mu_1, \mu_2 < 0$. Then, it follows from (4.1) and the Pohozaev identity that

$$P(u, v) = [u]_{H^s}^2 + [v]_{H^s}^2 - \|u\|_{2_s^*}^{2_s^*} - \|v\|_{2_s^*}^{2_s^*} - \eta_1\|u\|_{2_s^*}^{2_s^*} - \eta_2\|v\|_{2_s^*}^{2_s^*}$$
$$- \gamma 2_s^* \int_{\mathbb{R}^N} |u|^\alpha |v|^\beta = 0.$$

Again by the definition of weak solution for the above system (4.1), we have

$$[u]_{H^s}^2 + [v]_{H^s}^2 + \mu_1\|u\|_2^2 + \mu_2\|v\|_2^2 = \|u\|_{2_s^*}^{2_s^*} + \|v\|_{2_s^*}^{2_s^*} + \eta_1\|u\|_{2_s^*}^{2_s^*} + \eta_2\|v\|_{2_s^*}^{2_s^*}$$
$$+ \gamma 2_s^* \int_{\mathbb{R}^N} |u|^\alpha |v|^\beta.$$

The proof is now complete. □

Thus, we obtain

$$\mu_1\|u\|_2^2 + \mu_2\|v\|_2^2 = \mu_1 m_1^2 + \mu_2 m_2^2 = 0.$$

This is clearly a contradiction. The proof of Lemma 4.1 is complete.

**Proof of Theorem 1.2** It follows from Lemma 4.1 that the proof of Theorem 1.2 is complete. □

**Author Contributions** JZ and VR wrote the main manuscript text and they approved the final version of this paper.

**Funding** The research of Vicenţiu D. Rădulescu was supported by a grant of the Romanian Ministry of Research, Innovation and Digitization, CNCS/CCCDI-UEFISCDI, project number PCE 137/2021, within PNCDI III.

**Availability of data and materials** All data generated or analysed during this study are included in this article.

## Declarations